\documentclass[12pt,a4paper]{article}
\usepackage{amsmath}
\usepackage{amsfonts}
\usepackage{amssymb}
\usepackage{graphicx}
\usepackage{graphics}
\usepackage{makeidx}
\usepackage{natbib}
\input xypic.tex
\input xy
\xyoption{all}
\xyoption{2cell}
\UseAllTwocells

\def\AQ{{\mathbb A}}
\def\A{{\cal A}}
\def\BQ{{\mathbb B}}
\def\B{{\cal B}}
\def\AQ{{\mathbb A}}
\def\A{{\cal A}}
\def\BQ{{\mathbb B}}
\def\B{{\cal B}}
\def\C{{\cal C}}
\def\D{{\cal D}}
\def\G{{\cal G}}
\def\I{{\cal I}}
\def\L{{\cal L}}
\def\X{{\cal X}}
\def\MM{\underline{M}}
\def\End{{\bf End}\:}
\def\CAT{{\bf CAT}}
\def\CA2{\mathbb{CAT}}
\def\Const{\Delta}
\def\Cons{{\rm Const}}
\def\po{\circ}
\def\id{{\rm id}}
\def\Id{{\rm Id}}
\def\ran{{\rm Ran}}
\def\Ran{{\underline{\ran}}}
\def\ob{{\rm Ob}}
\def\mor{{\rm Mor}}
\def\imp{\Rightarrow}
\def\iff{\Leftrightarrow}
\def\hom{{\rm hom}}
\def\Hom{{\rm Hom}}

\def\barr{\begin{array}{rcl}}
\def\earr{\end{array}}

\def\alg{\!-\!{\bf alg}}
\def\Alg{{\bf Alg}\:}

\def\Con{{\bf Con}\:}

\newtheorem{theo}{Theorem}[section]
\newtheorem{coro}[theo]{Corollary}
\newtheorem{lm}[theo]{Lemma}

\newtheorem{note}[theo]{Note}
\newtheorem{rem}[theo]{Remark}
\newtheorem{de}[theo]{Definition}
\newtheorem{examp}[theo]{Example}
\newenvironment{proo}{{\bf {Proof:}}}{\hfill $\square$}

\begin{document}
\begin{center}
\LARGE{\bf Kan Extensions in Context of Concreteness }\\
\vspace{1 cm}
\normalsize{Jan Pavl\'ik} \\
\vspace{0.5 cm}
\small Faculty of Mechanical Engineering, Brno University of Technology,\\
Brno, Czech Republic\\
 {\tt pavlik@fme.vutbr.cz}\\
\end{center}

\begin{abstract} This paper contains results from two areas -- formal theory of Kan extensions and concrete categories. The contribution to the former topic is based on the extension of the concept of Kan extension to the cones and we prove that limiting cones create Kan extensions. The latter topic focuses on two significant families of concrete categories over an arbitrary category. Beck categories are defined by preservance properties while newly introduced l-algebraic categories are described by limits of categories of functor algebras. The latter family is shown to be rather natural.

The well known Beck's theorem states that the monadic categories are precisely the Beck categories with free objects. We strengthen this theorem by weakening the assumptions of the existence of free objects and we replace it by existence of some Kan extensions, namely the pointwise codensity monads. Moreover, using the result on Kan extensions of cones we show that for l-algebraic categories even weaker assumption fits.
\end{abstract}


Beck categories \citep[see][]{manes,Rac} are the concrete categories over some base category $\C$ such that their forgetful functor creates limits and certain colimits. Recall that the monadic categories are categories concretely isomorphic to Eilenberg-Moore categories for some monad over $\C$. Monadic categories are well-behaved categories, which tend to be categorical formulations of varieties of algebras.
The famous Beck's theorem states that monadic categories are precisely Beck categories with free objects, i.e., those Beck categories whose forgetful functor is right adjoint.

In this paper we show stronger versions of Beck's theorem, previously presented only in \citep{pd}. The main results involve notions of codensity monads, the concept of which is more general than that of adjunction. To explain the necessary background for codensity monads we recall the theory of Kan extensions. It is useful to derive a concept of Kan extension for the cones, which is only an instance of a more general concept definable in 2-categories. We prove a formal 2-categorical statement which implies that limiting cones create Kan extensions and we use this property in context of Beck categories and codensity monads.
There are two levels of Kan extensions (general and pointwise) which distinguishes two levels of needed conditions for Beck categories to be monadic. In order to express these conditions comprehensively, we introduce the notion of l-algebraic category. This is proved to be a very natural concept and a family of such concrete categories contains most of the categories with algebra-like objects.

Finally, we derive two characterizations of monadic categories which are stronger then the Beck's theorem. Namely, we show that monadic categories can be characterized precisely as
\begin{itemize}
\item Beck categories with pointwise codensity monad,
\item L-algebraic categories with codensity monad.
\end{itemize}

We work in Von Neumann -- Bernays -- G\"{o}del set theory with axiom of choice for classes and we use several levels of categories. The most usual case is when objects form a class and the morphisms between two objects form a set (i.e. the local smallness). Together with functors they form a category $\CAT$. Any structure of this or larger size will be referred to as a meta-structure. Generally, however, we use the prefix meta- only if we need to emphasize the actual size of the feature.

By ${\bf 1}$ we denote the category with only one element $0$ and identity while category ${\bf 2}$ contains two objects $0$ and $1$ and morphisms $\iota:0 \to 1$, $\id_0$, and $\id_1$.
\section{Kan Extensions and 2-categories}
To express the main results we need to work with the notions of universality and Kan extensions. Recall that, given a functor $R : \C \to \D$, an object $A \in \ob\D$ has an $R$-universal arrow with the base object $B \in \ob\C$,  iff
\begin{eqnarray}
\hom_{\D}(A,-)\po R \cong \hom_{\C}(B,-).\label{uni}
\end{eqnarray}
\begin{rem}\label{colim}
For instance, if $\D$ is a small category and $\Cons : \C \to \C^{\D}$ sends objects and morphisms on constant functors and constant transformations, respectively, then, given a diagram $D: \D \to \C$, its colimit (if any) is the base object of a $\Cons$-universal arrow over $D$.
\end{rem}
\subsubsection*{Kan Extensions and Codensity Monds}
Let us recall the concept of Kan extensions briefly. For the main sources, see \citep{mac,bor1}.

Let $\A$, $\B$ and $\C$ be categories and $S: \A \to \B$ , $U:\A \to \C$ be functors. A {\em right Kan extension of $S$ along $U$}
is a pair $\Ran_U S=(T,e)$ consisting of a functor $T: \C \to \B$ and a natural transformation $e:TU \to S$ satisfying the following universal property (Kan universality):\\
given a functor $T': \C \to \B$ and a transformation $e':T'U \to S$, then there is a unique transformation $t:T' \to T$ such that $e'=e \po tS$. Then we write $T=\ran_U S$.

$\Ran_U S$ can be defined equivalently as a $[-\po U]$-couniversal arrow $(T,e)$ over a functor $S: \A \to \B$, here $[-\po U]$ is the abbreviation of the corresponding functor $\CAT(U,\D): \CAT(\C,\D) \to \CAT(\A,\D)$. Hence the right Kan extension can be given by an isomorphism
\begin{eqnarray}
\hom_{\End \C}(-,S) \po [-\po U] &\cong& \hom_{\End \C}(-,\ran_U S). \label{homhom}\label{kan2}
\end{eqnarray}

The right Kan extension is said to be {\em pointwise} if it can be obtained by the following procedure:

Given an object $A$ in $\C$, then $Q_A:A \downarrow U \to \A$ will denote the forgetful functor for the comma category. Consider the functors $S \po Q_A$ for every $A \in \ob\C$. If each of these functors, seen as a large diagram, has a limit $T_A \in \ob\B$, then the assignment $A \mapsto T_A$ yields a functor $T:\C \to \B$ and there is a natural transformation $e:T\po U \to S$ such that $(T,e)=\Ran_U S$.

Suppose $\Ran_U S=(T,e)$ exists. Then is pointwise iff it is preserved by $\hom(C,-)$ for every object $C$ in $\C$.

The relation between Kan extensions and adjunction is captured in the equivalence of the following statements:
\begin{enumerate}
\item \label{adran1}$F$ has a left adjoint.
\item \label{adran2}$\Ran_F \Id_{\A}$ exists and $F$ preserves it.
\item \label{adran3} $\Ran_F \Id_{\A}$ exists and $F$ preserves all right Kan extensions with the values in $\A $.
\end{enumerate}

Since, on a category with copowers, every hom-functor is right adjoint, the implication $(\ref{adran1}) \imp (\ref{adran3})$, together with the previous statement, yields an important fact:
A right Kan extension with the values in a category with copowers is pointwise.

Given a diagram $D:\D \to \A$ and a terminal functor $T:\D \to {\bf 1}$, then $\Ran_T D=(C_L,\lambda)$ consists of the constant functor $C_L:{\bf 1} \to \A$
with $L$ being the limit object of $D$ and $\lambda=\{\lambda_d:L \to Dd|,d \in \ob\D\}$ being the limiting cone. Dually, colimits can be obtained by left Kan extensions.

Given a functor $U:\A \to \C$ with the right Kan extension $\Ran_U U=(M,e)$, then there exists a monad $\MM=(M,\eta,\mu)$ on $\C$ with the transformations $\eta:\Id_{\C} \to M$ and $\mu:M^2 \to M$ induced, via Kan universality, by $\id_U:U \to U$ and $e \po Me:M^2 U\to U$, respectively.

The above monad $\MM=(M,\eta,\mu)$ together with the transformation $e$ is called {\em codensity monad} for $U$. If the corresponding right Kan extension is
pointwise, we say $(\MM,e)$ is a {\em pointwise codensity monad}. If $(\A,U)$ is the concrete category, we say $\MM$ is its codensity monad.

\subsubsection*{Kan Extensions in 2-categories}\label{2-cat}
Given a 2-category $\AQ$, then, for each pair of its 0-cells $A$, $B$ the {\em category} of 1-cells $A \to B$ as objects and 2-cells as morphisms will be denoted by $\Hom(A,B)$ while we keep the notation of {\em hom-sets} (of 1-cells or 2-cells) denoted by "$\hom$".
For a more detailed overview of a 2-category theory, see \citep{street,lackstr,lack}.

The concept of Kan extensions is well established in every 2-category -- just think of every 0-cell, 1-cell, 2-cell as of a category, a functor, a natural transformation, respectively.

We are going to prove a property on Kan extensions in a 2-category which might be well-known but no source known to the author seems to contain it. It involves the notion of 2-couniversality (a weaker condition to an adjunction between 2-categories).

\begin{de}({\em 2-couniversality})
Given two 2-categories $\AQ$, $\BQ$, a 2-functor $L: \AQ \to \BQ$ and an object $B$ of $\BQ$, we say an object $Q$ of $\AQ$ is {\em $L$-2-couniversal} (briefly $L$-couniversal) over $B$ if there exists a 2-natural isomorphism between 2-functors $$K:\Hom_{\BQ}(-,B) \po L \to \Hom_{\AQ}(-,Q).$$
\end{de}
From now on, let the above instance of 2-couniversality take place.
\begin{lm}
Given objects $A$, $C$ of $\AQ$, morphisms $g:A \to C$ and $h:LC \to B$, then
$$K_A(h \po Lg)=K_Ch \po g.$$
\end{lm}
\begin{proo}
Since $K$ is a 2-natural transformation, we have $$K_A(h\po Lg)=K_A \po \Hom(Lg,B)(h)=\Hom(g,Q) \po K_C(h)=K_Ch \po g.$$
\end{proo}

\begin{rem}
Given objects $A$, $C$ of $\AQ$, morphisms $f:LA \to B$, $g:A \to C$ and $h:LC \to B$, then the 2-cell-function of $K_A$ defines an isomorphism
$\hom(h \po Lg,f) \cong \hom(K_C(h) \po g,K_A(f)) $
and due to naturality of $K$ we get a natural isomorphism
\begin{eqnarray}
\hom(-,f)\po [-\po Lg] &\cong& \hom(-,K_A f)\po [-\po g] \po K_C. \label{ka}
\end{eqnarray}

\end{rem}

\begin{theo}\label{kanuni}
Let $\AQ$ and $\BQ$ be 2-categories, $L:\AQ \to \BQ$ be a 2-functor and $B$ an object in $\BQ$ with an $L$-couniversal object $Q \in \ob\BQ$. Given objects $A, C$ in $\AQ$ and morphisms $f:LA \to B$, $g:A \to C$, then
$$\Ran_{g}(K_A f) \cong K_C \Ran_{Lg} f$$
whenever $\Ran_{Lg} f$ exists.
\end{theo}
\begin{proo}
Let $\ran_{Lg} f$ exist, then by application of (\ref{kan2}) and of the previous remark
we have:
$$
\barr
\hom(-,K_A f)\po [-\po g] \po K_C &\stackrel{(\ref{ka})}{\cong}&  \hom(-,f)\po [-\po Lg] \\
&\stackrel{(\ref{kan2})}{\cong}&  \hom(-,\ran_{Lg}f)\\
&\stackrel{(\ref{ka})}{\cong}&  \hom(-,K_C \ran_{Lg}f) \po K_C
\end{array}
$$
and since $K_C$ is an isomorphism, we get $$\hom(-,K_A f)\po [-\po g] \cong \hom(-,K_A \ran_{Lg}f).$$ Hence $K_A \ran_{Lg}f \cong \ran_{g}(K_A f)$ which can be extended analogously on natural transformations to get the required equality.
\end{proo}

\subsubsection*{Kan Extensions of Cones}
Our interest will be in limits of concrete categories. Therefore, it will be suitable to study Kan extensions of cones. These will be defined using the language of 2-categories.

Consider a 2-category $\AQ$. If $\D$ is a small category, then the category $\AQ^{\D}$ of functors $\D \to \AQ$ can be treated again as a 2-category with the cell-structure defined pointwise. Now the objects of $\AQ^{\D}$ can be seen as $\D$-domained diagrams.

We define a 2-functor $\Const: \AQ \to \AQ^{\D}$ as depicted bellow (with 0-cells $\A$, $\B$, 1-cells $F$, $G$ and a 2-cell $\alpha$):

$$ \begin{array}{ccc}
\xymatrix{
\A \rrtwocell_G^F{_\;\:\alpha} &&
\B
} & \mapsto &
\xymatrix{
C_{\A} \rrtwocell_{c_{G}}^{c_{F}}{_\quad{\rm c}_{\alpha}} &&
C_{\B}
}
\end{array}
$$
where $C_{\A}$ denotes an $\A$-constant functor, $c_F$ is a $F$-constant natural transformation and ${\rm c}_{\alpha}$ stands for an $\alpha$-constant modification of natural transformations.

Given a category $\A$, then a 1-cell $\Const\A \to D$ is a $D$-compatible cone and a 2-cell between such cones is a $\D$-indexed collection of natural transformations.

By analogy to the Remark \ref{colim}, we can define cone limits in 2-categories using the notion of 2-couniversality.

\begin{de}
Given a 2-category $\AQ$, we say a diagram $D:\D \to \AQ$ has a limit $l \in \ob\AQ$ if $l$ is a $\Const$-2-couniversal object over $D$.
\end{de}

The concept holds obviously even for large 2-categorical structures. Consider a 2-metacategory $\CA2$ with the usual cell-structure and a category $\D$ (the smallness condition will not be necessary).
Now, given a diagram $D:\D \to \CAT$, functor $G:\A \to \B$ and a cone $F:\Const\A \to D$, a {\em Kan extension of the cone $F$ along the functor $G$} is the Kan extension of $F$ along $\Const G$ in the 2-category $\CA2^{\D}$.
$$\xymatrix{\Const\A \ar[rr]^{F}\ar[drr]^{\Const G} && D\\ && \Const\B \ar[u]_{\ran_{\Const G}F}}$$
Now the Theorem \ref{kanuni} has the following consequence:
\begin{lm}\label{konlim}
Let $D:\D \to \CA2$ be a diagram and $\L$ be its limit and $\A$ be a category. Given a cone $F:\Const\A \to D$ and a functor $G:\A \to \B$, then if $\ran_{\Const G} F$ exists, then
$$\ran_{G} \widetilde{F} \cong \widetilde{\ran_{\Const G} F}$$
where $\widetilde{\:}$ stands for the factorization over the limit cone.
\end{lm}
\begin{proo}
Since $\L$ is the limit of $D:\D \to \CA2$, it is a $\Const$-2-couniversal object over $D$, hence the statement follows directly form the Theorem \ref{kanuni}.
\end{proo}

We say $\Ran_{\Const G}F$ is {\em componentwise} Kan extension if it is a collection of the right Kan extensions, i.e., $$(\forall d\in \ob\D)\: (\Ran_{\Const G}F)_d= \Ran_{G} F_d.$$
\begin{rem}\label{compkancone}
Under the above assumptions, it is easy to show the equivalence of the following:
\begin{enumerate}
\item A componentwise right Kan extension of the cone $F$ along $G$ exists.
\item For every object $d \in \ob\D$, a $\Ran_G D_d$ exists and these Kan extensions are preserved by all composable functors of the form $D\phi$, $\phi \in \mor\D$.
\end{enumerate}
\end{rem}
Since the usual limits can be obtained as right Kan extensions of the diagram along the terminal functor, we may consider the following analogue.
\begin{rem}
Let there be a diagram $D:\D \to \CAT$ with the limit $\L$ and the limiting cone $L$ and a diagram $G:\G \to \L$ and the terminal functor $T:\G \to \1$. Coherently with the usual concept of Kan extensions, the $\Ran_{\Const T} (L \po \Const G)$ will be called a {\em limit of the cone diagram $L \po \Const G$}. In fact, it is a choice of the limit (object and the limiting cone) of each $L_d \po G$, such that they are preserved by functors $D\phi:Dd \to De$ for $\phi:d\to e$ in $\mor\D$.
If $D$ has a limit, any cone can be seen as a cone diagram since it factorizes uniquely over the limiting cone.
\end{rem}

As a direct consequence of the Lemma \ref{konlim} we can show that Kan extensions are created by limiting cones of categories.
\begin{theo}\label{limkancone}{\rm (Limiting cones create Kan extensions.)}
Let there be a diagram $D:\D \to \CAT$ with the limit $\L$ and the limiting cone $L=\{L_d:\L \to Dd|,d \in \ob\D\}$. Let there be categories $\A$, $\B$ and functors $V:\A \to \B$, $S:\A \to \L$. Consider the induced cone $K=L \po \Const S$. If $\Ran_{\Const V} K$ exists, then $\Ran_V S$ exists as well and $L \po \Const\Ran_V S \cong \Ran_{\Const V} K$.
\end{theo}
\begin{rem}\label{leftkancone}
By duality of categories $\A,\B,\L,D_d$ in the theorem, one gets that limiting cones create left Kan extensions as well.
\end{rem}
\begin{lm}\label{limkan}
Let there be a diagram $D:\D \to \CAT$ with the limit $\L$ and the limiting cone $L$ and $G:\G \to \L$ be a diagram in $\L$. If the limit of the cone diagram $L \po \Const G$ exists, the limit of $G$ exists as well and $L$ preserves it.
\end{lm}
\begin{proo}
Observe that, for every pair of functors $\xymatrix{\bullet & \bullet \ar[l]_{Q} \ar[r]^{P}& \bullet}$ there is an isomorphism $\Ran_{\Const Q}\Const P \cong \Const \Ran_Q P$ and that the limit of the constant functor is the target object itself.
Then, since $\widetilde{H}=L \po \Const H$ for every functor $H:\X \to \L$,
we have
$$\barr
\lim(L \po \Const G) &=&\ran_{\Const T}(L \po \Const G) =\ran_{\Const T}(\widetilde{G})\\
&\cong& \widetilde{\ran_T G}\\
&=&L \po \Const \ran_T G \cong L \po \ran_{\Const T} \Const G\\
&=&L \po \lim \Const G = L \po \Const \lim G.
\earr
$$
\end{proo}
\section{Concrete Categories and Algebras}\label{limcc}
The reader is expected to be familiar with the concept of concreteness, functor algebras and monadicity. To clarify the notation, given a category $\C$, a $\C$-concrete category with the forgetful functor $U:\D \to \C$ will be denoted by $(\D,U)$ or just by $\D$ if we do not need to emphasize the name of the forgetful functor. If the choice of this functor is obvious, the forgetful functor is usually denoted by $U_{\D}$. If concrete categories $\A$ and $\B$ are concretely isomorphic, we write $\A \cong_{\C} \B$.

Recall from \citep{Rac}, a $\C$-concrete category $(\A,U)$ is said to be a {\em Beck category} if $U$ creates all limits and coequalizers of $U$-absolute pairs. Beck categories include all categories of algebras and of algebra-like objects.

Given a functor $F: \C \to \C$, the category of $F$-algebras will be denoted by $\Alg F$. A category isomorphic to $\Alg F$, for some functor $F:\C \to \C$, will be called {\em f-algebraic}. It can be shown that every f-algebraic category is Beck.

\begin{rem}\label{alglim}
Given a category $\C$, the metaclass of $\C$-concrete categories and $\C$-concrete functors will be denoted by $\Con \C$ and referred to as {\em category of $\C$-concrete categories}.
The operator $\Alg$ may be considered contravariant functor $\End \C \to \Con \C$  defined on the category of endofunctors on $\C$. It assigns, to a natural transformation $\phi:G \to F$, the concrete functor $\Alg\phi:\Alg F \to \Alg G$ given by equality $(\Alg\phi)(A,\alpha)=(A,\alpha \po \phi_A)$.
\end{rem}

To revise the notion of monadicity, recall that, given a monad $\MM= (M, \eta, \mu)$ on $\C$, a {\em category $\MM\alg$ of
$\MM$-algebras} (also called Eilenberg-Moore category) is a full subcategory of $\Alg M$ consisting of all $M$-algebras $(A,\alpha)$ satisfying the {\em
Eilenberg-Moore identities} $\alpha \po \eta_A=\id_A$,
$\alpha \po M\alpha =\alpha \po \mu_A$.
A category concretely isomorphic to $\MM\alg$ for some monad $\MM$ is called {\em monadic}. The Eilenberg-Moore category $\MM\alg$ has free objects and the corresponding free functor $W_{\MM\alg}$ assigns, to an object $A$, the algebra $(M A,\mu_A)$.
It yields the adjunction $W_{\MM\alg} \dashv U_{\MM\alg}$ with the associated monad equal to $\MM$, which is also the codensity monad for $U_{\MM\alg}$. For more detailed treatment of monadic categories see \citep{mac,AHS,bor1,bor2}.

The following well known theorem characterizes all monadic categories:
\begin{theo}{\rm (Beck's theorem)}\\
Let $\A$ be a concrete category. Then
\begin{center}
$\A$ is monadic $\iff$ $\A$ is a Beck category with free objects.
\end{center}
\end{theo}

\subsection{Limits of Concrete Categories}

In this section we will work with diagrams of concrete categories over a base category $\C$.
By a (small) {\em diagram $D$ of concrete categories} we mean a functor $D: \D \to \Con \C$ from some (small) category $\D$.

Since, for each concrete category $(\A,U)$, $U$ is a faithful functor $\A \to \C$, we may consider $\Con \C$ a subcategory of the "slice-category" $\CAT / \C$ containing those pairs $(\A,U)$ where $U$ is faithful.

\begin{lm}
$\Con \C$ is closed under limits in $\CAT / \C$.
\end{lm}
\begin{proo}
The limits of concrete categories are computed as limits in a slice category -- the product is a wide pullback of the diagram connected by forgetful functors and the equalizer is an equalizer of the "non-concrete part" of the pair. The forgetful functor $U_{\L}$ for the limit category $\L$ can be obtained by composition of any forgetful functor of a category in a diagram with the corresponding limit-cone component. All we need to show is the faithfulness of this forgetful functor.

Let a pair $\xymatrix{\bullet \ar@<+0.5ex>[r]\ar@<-0.5ex>[r] &\bullet}$ in $\L$ be collapsed by $U_{\L}$ into an arrow. Then, due to factorization of $U_{\L}$ over faithful functors for each component of the corresponding diagram we have a cone with the domain ${\bf 2}$. Hence, there is a unique functor ${\bf 2} \to \L$ which makes the original arrows equal.
\end{proo}

\begin{rem}\label{conpro}
Since $\CAT$ is complete w.r.t. small diagrams, so is $\CAT / \C$, hence $\Con \C$ is complete as well. In particular, there exist products in $\Con \C$. Given two concrete categories $\A$, $\B$, their concrete product will be denoted by $\A \times_{\C} \B$, and can be obtained as a fibre-wise product, i.e., $$U_{\A \times_{\C} \B}^{-1} (A)=U_{\A}^{-1} (A) \times U_{\B}^{-1}(A).$$ The product of an infinite number of concrete categories can be described analogously.
The terminal object, i.e., the product over the empty index set, is the base category itself and the terminal morphisms are forgetful functors.
\end{rem}
It is easy to check that the functor $\Alg$ turns colimits into limits, i.e., e.g., $$\Alg(F+G) \cong \Alg F \times_{\C} \Alg G.$$
\begin{rem}\label{classint}
It is easy to see, that, moreover, the intersection of a large collection of subcategories exists.
Clearly, it contains all objects and morphisms, which occur in all categories of collection.
\end{rem}

\begin{lm}\label{limbeck}
Let $(\L,U)$ be the concrete limit of a diagram $D:\D \to \Con \C$. Let $G:\G \to \L$ be a diagram in $\L$ such that the forgetful functor $U_d$ creates the limit of $L_d \po G$ for each $d \in \ob\D$. Then $U$ creates the limit of $G$.

Dually, if the forgetful functors create the corresponding colimits, then so does $U$.
\end{lm}
\begin{proo}
To see this property for limits, let $U \po C_G$ have a limit. Since $U=U_d \po L_d$ for every object $d \in \D$ and the functors $U_d$ create the limit of $L_d \po G$, we have the limit of each $L_d \po G$. Moreover, the functors $L_d$ create the limiting cones, i.e., there is always a unique cone $\lambda^d_q$ of $L_d G$ for each $d$. Hence, for every morphism $\phi:d \to e$ in $\D$, we have $D\phi \lambda^d_q=\lambda^e_q$. Therefore there is a $D$-compatible cone $M$ with the domain $\1$ such that $M_d(0)=\lim L_d G$. Hence $M$ is a limit of the cone diagram $L \po C_G$. Due to Lemma \ref{limkan} it creates the limit of $G$, i.e., the limit of $G$ exists and is preserved by every $L_d$.

The property for colimits can be reached using the dual property (Remark \ref{leftkancone}).
\end{proo}

Immediately, we get the consequence.
\begin{coro}\label{limbeck}
A (possibly large) concrete limit of Beck categories is Beck.
\end{coro}
\subsection{Algebras for a Concrete Diagram}
\begin{de}
A limit, denoted by $\Alg D$, of a (possibly large) diagram $D: \D \to \Con \C$ with every object mapped on an f-algebraic category, will be called an {\em l-algebraic category}. If the category $\D$ has a weakly initial object, we say $\Alg D$ is {\em homogenous}.
\end{de}
\begin{rem}
Let $D :\D \to \End \:\C$ be a diagram with $D(x)=\Alg F_x$, where $F_x:\C \to \C$ is a functor for every object $x \in \D$.
The objects of category $\Alg D$ may be seen as {\em algebras for the diagram $D$}, i.e., the collections of
$\C$-morphisms $\{\alpha_x:F_xA \to A|x \in \D\}$ satisfying, for each $f:x \to y$ in $\D$, the $D$-compatibility condition $D(f)(A,\alpha_x)=(A,\alpha_y)$. The morphisms in $\Alg D$ are the morphisms of algebras for each
$x$, i.e., $\phi:(A,\alpha)\to (B, \beta)$ is a morphism if $\phi \circ \alpha_x = \beta_x \circ
F_x(\phi)$ for every $x \in \ob\D$.
\end{rem}
An analogous approach to the algebras for a {\em diagram of monads} was studied by Kelly in \citep{Kelly}.

The concrete product $\Alg F \times_{\C} \Alg G$ for some endofunctors $F, G$ on $\C$ is an example of l-algebraic category. Its objects are $\C$-objects with two structure arrows -- for both $F$ and $G$. However, if $\C$ has coproducts, then the category $\Alg F \times_{\C} \Alg G$ is isomorphic to $\Alg (F + G)$. Generally, if $\C$ is cocomplete and the diagram in the definition maps all morphisms on "f-algebraic functors", i.e., those concrete functors which are $\Alg$-images of natural transformations, then the obtained category is an f-algebraic category. Therefore we will focus especially on the diagrams which do not factorize over functor $\Alg$.

Here we show the basic property of each l-algebraic category.

\begin{lm}\label{lbeck}
Every l-algebraic category is Beck.
\end{lm}
\begin{proo}
The statement is a direct consequence of property of f-algebraic categories and Lemma \ref{limbeck}.
\end{proo}

In Example \ref{proti} we show a Beck category which is not l-algebraic. Hence, the above implication cannot be reversed.

In \citep{pd}, the author shows that if $\C$ is cocomplete, all varieties in sense of \citep{AP} and all algebraic categories in sense of \citep{Rac} are l-algebraic.
\subsection{Polymeric Varieties}\label{polysect}
Another important family of concrete categories was introduced in \citep{p1}. It was shown that it includes many natural examples. To recall the definition, we need the notion of {\em polymer}.
\begin{de}
Let $(A,\alpha)$ be an $F $-algebra. Given $n \in \omega$, a {\em $n$-polymer} of an algebra $(A,\alpha)$
is the morphism $\alpha^{(n)}:F^n(A) \to A$ in $\C$ defined recursively:
$$\alpha^{(0)}=\id_A, \alpha^{(n+1)}=\alpha \circ F \alpha^{(n)}.$$
\end{de}
\begin{rem}\label{polyfun}
The assignment $(A,\alpha) \mapsto (A,\alpha^{(n)})$ clearly defines a functor $P_n:\Alg F \to \Alg F^{n}$.
\end{rem}
Now we have a sequence of functors $F^n$ together with functors between the categories of their algebras. This leads us to the definition:
\begin{de}\label{polymeric}
Let $n \in \omega$ and $G $ be an endofunctor on $\C$. A natural transformation $\phi:G  \to F^n$ is called
{\em $n$-ary polymeric $G$-term} in
category of $F$-algebras. A pair $(\phi,\psi)_p$ of polymeric $G$-terms of arities $m$,$n$,
respectively, is called {\em polymeric identity} of {\em
arity-pair} $(m,n)$ with {\em domain} $G$.
Moreover, for an $F$-algebra $(A, \alpha)$, we define $$(A, \alpha)\models (\phi,\psi)_p \quad
\stackrel{def}{\iff} \quad \alpha^{(m)} \circ \phi_A = \alpha^{(n)} \circ \psi_A,$$ and we say that the $F$-algebra $(A, \alpha)$ {\em satisfies} the
polymeric $G$-identity $(\phi,\psi)_p$.

For a class $\I$ of polymeric identities we define a {\em polymeric variety of $F$-algebras} as the class of all algebras
satisfying all $(\phi,\psi)_p \in \I$. The corresponding full subcategory of $\Alg F$ is denoted
by $\Alg (F,\I)$. If $\I$ is a singleton, we say that $\Alg (F,\I)$ is {\em single-induced}.
A category concretely isomorphic to $\Alg (F,\I)$ for some $F$ and $\I$ will be called {\em polymeric}.
\end{de}
\begin{lm}
Every polymeric category is homogenous l-algebraic category.
\end{lm}
\begin{proo}
It is easy to see that a single-induced polymeric variety $\Alg (F,(\phi,\psi)_p)$, where $(\phi,\psi)_p$ is $(m,n)$-ary polymeric $G$-identity, is an equalizer of the pair of concrete functors obtained by the following compositions:
$$\xymatrix{&\Alg F^m \ar[dr]^{\Alg \phi}&\\
\Alg F \ar[ur]^{P_m}\ar[dr]^{P_n} &&\Alg G\\
&\Alg F^n. \ar[ur]^{\Alg \psi}&\\
}$$
A polymeric variety induced by a class of polymeric identities is the intersection of polymeric varieties induced by single polymeric identities, hence a homogenous l-algebraic category.
\end{proo}

The following lemma shows the presentation of monadic categories by polymeric identities.
\begin{lm}\label{monad}
\label{poly-monad}
Every monadic category is l-algebraic.
\end{lm}
\begin{proo}
Given a monad $\MM=(M ,\eta,\mu)$, the Eilenberg-Moore category $\MM\alg$ is a
polymeric variety of $M $-algebras induced by polymeric identities $(\eta, \id_{\Id})_p,
(\id_{\MM^2},\mu)_p$ of domains $\Id$, $\MM^2$, respectively, and arity-pairs $(1,0), (2,1)$,
respectively.
$$\xymatrix{
M^1&\Id \ar[l]_{\eta} \ar[r]^{\id}& M^0,
}\quad \xymatrix{
M^2&M^2 \ar[l]_{\id_{M^2}} \ar[r]^{\mu} &M^1.} $$
Therefore $\MM\alg$ is a polymeric category, hence it is l-algebraic.
\end{proo}
\section{Concreteness and Universality}
In order to solve some problems involving categories of algebras on category $\C$ such
as existence of free objects, it might be useful to know their properties on the level of
objects of category $\Con \C$. Their status can be expressed in terms of properties of contravariant functor $\Alg$.
To simplify the proofs, we express the concept of universality in the language of Kan extensions and codensity monads. The main results characterize free-objects existence in l-algebraic categories and all Beck categories and provide another characterization of monadic categories.

The author acknowledges the advice from H. E. Porst, J. Rosick\'{y} and J. Velebil concerning the connection of universality and the concept of Kan extensions.

All the topic in this chapter deals with the category $\C$, generally without any additional assumptions.
\subsection{$\Alg$-universality and Kan Extensions}
\begin{rem}\label{transalg}
Let $F$ be an endofunctor on $\C$, $U_F:\Alg F \to \C$ be the forgetful functor and $\A$ be a category. Then each functor $H:\A \to \Alg F$ can be seen as $H=(K,\kappa)$ where $K=U_F \po H: \A \to \C$ and $\kappa:FK \to K$ is a natural transformation such that $HA=(KA,\kappa_A)$ for every object $A$ in $\A$. Moreover, if $(\A,U_{\A})$ is a $\C$-concrete category and $H$ is concrete, then $K=U_{\A}$.
This property can be easily expressed in terms of isomorphism between two (illegitimate) contravariant functors
\begin{eqnarray}
\hom_{\Con \C}(\A,-)\po \Alg &\cong& \hom_{\End \C}(-,U_{\A}) \po [-\po U_{\A}], \label{funhom1}
\end{eqnarray}
where $[-\po U_{\A}]=\hom_{\CAT}(U_{\A},\C)$ is given by the assignment $F \mapsto F \po U_{\A}$ (see the introduction to Kan extensons).
\end{rem}
We will apply the concept of universal arrows on the contravariant functor $\Alg: \End \C \to \Con \C$ as follows:\\
Given a $\C$-concrete category $\A$ and a $\C$-endofunctor $F$ with a concrete functor $H:\A \to \Alg F$, then $(F,H)$ is an $\Alg$-universal arrow for $\A$ iff for every $\C$-endofunctor $G$ with a concrete functor $J: \A \to \Alg G$ there exists a unique transformation $\widetilde{J}:G \to F$ such that $J= \Alg \widetilde{J} \po H$. This can be expressed by the isomorphism
\begin{eqnarray}
\hom_{\Con \C}(\A,-) \po \Alg  &\cong& \hom_{\End \C}(-,F). \label{funhom2}
\end{eqnarray}
The following observation is due to J. Rosick\'{y}:
\begin{lm}
Let $\A$ be a concrete category. Then
\begin{center}
$\A$ has an $\Alg$-universal arrow $\iff$ $\A$ has a codensity monad.
\end{center}
\end{lm}
\begin{proo}
Since the existence of a codensity monad is just a different way of saying that $\Ran_{U_{\A}} U_{\A}$ exists, we need to show that $\Alg$-universal arrow exists for $\A$ iff $\Ran_{U_{\A}} U_{\A}$ does. But the natural isomorphisms
(\ref{funhom1}) and (\ref{funhom2}) together with (\ref{homhom}) yield exactly what we need: a functor $F$ is a base of an $\Alg$-universal arrow for $\A$ iff $F=\ran_{U_{\A}} U_{\A}$ since
$$\hom(-,F)\stackrel{(\ref{funhom2})}{\cong}  \hom_{\Con \C}(\A,-) \po \Alg \stackrel{(\ref{funhom1})}{\cong} $$
$$\stackrel{(\ref{funhom1})}{\cong}\hom_{\End \C}(-,U_{\A}) \po \hom_{\CAT}(U_{\A},\C)
\stackrel{(\ref{homhom})}{\cong} \hom(-,F)$$
\end{proo}

This fact enables us to investigate the universality in terms of Kan extensions and codensity monads. In the following sections, we will use this language to characterize the free-objects existence for two significant families of concrete categories.
\subsection{Beck Categories with Codensity Monad}
\subsubsection{Categories with Pointwise Codensity Monad}
\begin{lm}
Let $U:\A \to \C$ be a functor which creates all limits and has a pointwise codensity monad. Then $U$ has a left adjoint.
\end{lm}
\begin{proo}
Let $(M,\eta)$ be the pointwise codensity monad for $U$. Then, for every $A \in \C$, $MA=\lim(U\po Q_A)$ where $Q$ is the projection functor $A\downarrow U \to \A$. Since $U$ creates the limits, there exists an object $L(A)=\lim Q_A$ in $\A$ such that $\lim(U\po Q_A)=U\po \lim Q_A$. Moreover, the whole limiting cone can be lifted to $\A$, i.e., for every morphism $f:A \to UB$ there is a unique morphism $\widetilde{f}:L(A) \to B$ such that $U\widetilde{f}$ is $\langle f \rangle$, i.e., the $f$-labeled component of the limiting cone for $U \po Q_A$. This implies the equality $f=U\widetilde{f} \po \eta_A$ for every $f$. Hence $\eta_A:A \to MA =UL(A)$ is the $U$-universal arrow for $A$. Therefore the assignment $A \mapsto L(A)$ can be extended to a functor $L:\C \to \A$, which is left adjoint to $U$.
\end{proo}

As a direct consequence we have the following.
\begin{coro}
Every Beck category with a pointwise codensity monad is monadic.
\end{coro}
This, in fact, implies a stronger version of Beck's theorem:
\begin{theo}
Let $\A$ be a concrete category. Then
\begin{center}
$\A$ is monadic $\iff$ $\A$ is a Beck category with a pointwise codensity monad.
\end{center}
\end{theo}
We emphasize that if the category has colimits, then each codensity monad is pointwise.
 Hence, in case such a case, for every Beck category, the existence of free objects is fully determined by the existence of a codensity monad.
 Particularly in case of varieties, there is another description of free objects derived in \citep{p}.
\subsubsection{L-algebraic Categories with Codensity Monad}
In the text bellow we show a seemingly weaker result, that every l-algebraic category with a codensity monad has free objects. However, this framework is independent of pointwiseness of a given codensity monad and cannot be derived from the above theorem, as shown at the end of this chapter.
\begin{lm}\label{ranconc}
Let $(\A,U)$ be a $\C$-concrete category with the right Kan extension of $U$ along itself. Then for every $\C$-endofunctor $F$ and a concrete functor $J:\A \to \Alg F$ there exists $\Ran_U J =(V,e)$. This Kan extension is preserved by
\begin{enumerate}
\item the forgetful functor $U_F:\Alg F \to \C$,
\item every concrete functor $T: \Alg F \to \Alg G$ (for every $G:\C \to \C$).
\end{enumerate}
\end{lm}
\begin{proo}
Let $J:\A \to \Alg F$ be a concrete functor and $\Ran_{U}U=(M,\epsilon)$. Due to the Remark \ref{transalg}, $J=(U,\iota)$ for some natural transformation $\iota:FU \to U$. Then $\epsilon:MU \to U$ yields the natural transformation $\iota \po F\epsilon:FMU \to U$. Hence there is a transformation $\upsilon:FM \to M$ such that the diagram
$$\xymatrix{
MU \ar[rr]^{\epsilon} &&U\\
FMU \ar[u]_{\upsilon U} \ar[rr]^{F\epsilon} && FU \ar[u]^{\iota}}$$
commutes.
We have gained a functor $V=(M,\upsilon):\C \to \Alg F$. The diagram above implies that the transformation $\epsilon$ can be extended to a transformation $\zeta :VU \to J$ such that $U_F \zeta=\epsilon$.

Let $G$ be another $\C$-endofunctor and let $U_G:\Alg G \to \C$ be the forgetful functor. Let $T=(U_F,\tau):\Alg F \to \Alg G$ be a concrete functor. We will show, that $(TV,T \zeta)$ is the right Kan extension of $TJ=(U,\tau J)$ along $U$.

Consider a functor $K=(U_G K, \kappa):\C \to \Alg G$ and a transformation $\lambda:KU \to TJ$. Then $U_G\lambda:U_G K U \to U$ is a natural transformation, hence there is $\rho:U_G K \to M$ such that $U_G \lambda = \epsilon \po \rho U$. Since $\lambda$ is the transformation $(U_G K U, \kappa U) \to (U_G TJ, \tau J)$, we have $U_G \lambda \po \kappa U = \tau J \po GU_G \lambda$.
The functor $T$ turns the diagram above into:
$$\xymatrix{
MU \ar[rr]^{\epsilon}&&U\\
GMU \ar[rr]^{G\epsilon}\ar[u]_{\tau VU} &&GU, \ar[u]^{\tau J} \\
}$$
hence $\tau J  \po G\epsilon = \epsilon \po \tau VU $.
We prove that $\rho$ underlies a transformation $K \to TV$.
$$\begin{array}{rcl}
\epsilon \po \rho U \po \kappa U &=& U_G \lambda \po \kappa U\\
&=& \tau J \po G U_G \lambda \\
&=& \tau J  \po G\epsilon \po G\rho U \\
&=& \epsilon \po \tau VU \po G\rho U.
\end{array}$$
Now from $(M, \epsilon)$ being the right Kan extension of $U$ along $U$, the uniqueness of factorization of a transformation $G U_G K U \to U$ over $\epsilon$ yields $\rho \po \kappa =\tau V \po G\rho$. Hence there is a transformation $\gamma:K \to TV$ such that $U_G \gamma = \rho$. Hence we have $U_G \lambda = \epsilon \po U_G \gamma U = U_G T\zeta \po U_G \gamma U$ which, due to faithfulness of $U_G$, implies
 $\lambda=T\zeta \po \gamma U$ and $(TV,T\zeta)$ becomes the right Kan extension of $TJ$ along $U$.

Since the procedure holds for every $G$ and $T$, by choice of $G=F$ and $T=\Id_{\Alg F}$ we get $(V,\zeta)=\Ran_U J$. For every other choice of $G,T$ we have $T\Ran_U J=T(V,\zeta)=(TV,T\zeta)=\Ran_U (TJ)$, hence $T$ preserves $\Ran_U J$.

Moreover $U_F\Ran_U J=U_F (V,\zeta)=(U_F V, U_F \zeta)=(M,\epsilon)=\Ran_U U$, hence even $U_F$ preserves $\Ran_U J$.
\end{proo}
\begin{note}
As pointed out by J. Velebil, since the whole proof is based only on the properties of functors and natural transformations, one can prove ana\-lo\-gous statement formally on the level of 2-categories.\end{note}
\begin{lm}\label{ranlim}
Let $(\A,U)$ be the limit of a diagram $D:\D \to \Con \C$, where $Dd$ is an f-algebraic category for each object $d$ of $\D$. If the codensity monad for $U$ exists, then $\Ran_U \Id_{\A}$ exists and is preserved by each component of the limiting cone $L_d:\A \to Dd$ $(d \in
\ob \D)$.
\end{lm}
\begin{proo}
If $\Ran_U U$ exists, then, for every $\D$-morphism $\phi:d \to d'$, the functor $D\phi:Dd \to Dd'$ is concrete between
f-algebraic categories, thus, due to Lemma \ref{ranconc}, $\Ran_U L_d$ and $\Ran_U L_{d'}$ exist and $\Ran_U
L_{d'}=D(\phi)\Ran_U L_d$. This is, due to the Remark \ref{compkancone}, equivalent to the existence of a componentwise right Kan extension of the cone $L$ along $U$. Then, by Theorem \ref{limkancone}, $\Ran_U \Id_{\A}$ exists and it is preserved by $L$.
\end{proo}
\begin{lm}
The forgetful functor $U$ of an l-algebraic category in $\Con \C$ with codensity monad creates the right Kan extension
of identity along $U$.
\end{lm}
\begin{proo}
If $(\A,U)$ is a limit of an empty diagram, then it is the terminal object of $\Con \C$, hence
$(\A,U)\cong (\C,\Id_C)$ and the situation is trivial.

Suppose $(\A,U)$ is a limit of a nonempty diagram and $\Ran_U U$ exists. The limiting cone is formed by f-algebraic categories $(\B_d,U_d)$ and concrete functors $L_d:\A \to \B_d$. Then according to the previous lemma, $\Ran_U Id$ exists and each functor $L_d$ preserves it. Moreover, due to Lemma \ref{ranconc}
each $U_d$ and the limiting cone $L_d$ preserves $\Ran_U L_d$,
hence their composition, which is equal to $U$, preserves $\Ran_U \Id$. Indeed, $U \Ran_U \Id =U_d L_d\Ran_U \Id \stackrel{\ref{ranlim}}{=}U_d
\Ran_U L_d \stackrel{\ref{ranconc}}{=} \Ran_U U_d L_d=\Ran_U U$.
\end{proo}

As a direct consequence of this lemma and the property of codensity monads we get the result:
\begin{theo}\label{codlcat}
Let $\A$ be an l-algebraic category. Then
\begin{center}
$\A$ has a codensity monad $\iff$ $\A$ has free objects.
\end{center}
\end{theo}

Since every monadic category is l-algebraic (see Lemma \ref{monad}), we get an alternative characterization of monadic categories:
\begin{theo}\label{monlcat}
Let $\A$ be a concrete category. Then
\begin{center}
$\A$ is monadic $\iff$ $\A$ is an l-algebraic category with a codensity monad.
\end{center}
\end{theo}
The statement cannot be extended to all Beck categories
unless the codensity monad is pointwise, as shown in the following example.
\begin{examp}\label{proti}
Consider the category $\C={\bf 2} + {\bf 2}$ consisting of objects $0,1,0',1'$ and morphisms $\iota:0 \to 1$, $\iota':0' \to 1'$ and identities. Let  $\A={\bf 1} + {\bf 1}$  and $U:\A \to \C$ be the inclusion of $\{0,0'\}$. Then the following holds:
\begin{enumerate}
\item \label{proti1}$(\A,U)$ is a Beck category.
\item \label{proti2}$(\A,U)$ has a codensity monad (the trivial monad).
\item \label{proti3}$U$ does not have an adjoint (1 does not have an universal arrow).
\end{enumerate}
\end{examp}

As a consequence, we see that there exists a codensity monad which is not pointwise and a Beck category that is not l-algebraic.
\section{Conclusion}
We have proved the propositions which, together with Beck's theorem, may be collected to the following theorem.
\begin{theo}
Let $\C$ be a category and $(\A,U)$ be concrete category over $\C$. The following statements are equivalent:
\begin{enumerate}
\item $\A$ is monadic.
\item $\A$ is Beck and $U$ has a left adjoint.
\item $\A$ 
is Beck and $U$ has a pointwise codensity monad.
\item $\A$ is l-algebraic and $U$ has a codensity monad.
\item $\A$ is an l-algebraic category with an $\Alg$-universal arrow.
\end{enumerate}
\end{theo}
\begin{rem}
If $\C$ has copowers, then we have even stronger result:
\begin{center}
$\A$ is a monadic category $\iff$ $\A$ is a Beck category with a codensity monad $\iff$ $\A$ is a Beck category with an $\Alg$-universal arrow.
\end{center}
\end{rem}
\paragraph{Acknowledgement} The author acknowledges partial support from Mi\-ni\-stry of Education of the Czech Republic, research plan no. MSM0021630518.

\begin{thebibliography}{}

\bibitem[Ad\'{a}mek et~al., 1990]{AHS}
Ad\'{a}mek, J., Herrlich, J.~H., and Strecker, G. (1990).
\newblock {\em Abstract and Concrete Categories: The Joy of Cats}.
\newblock John Wiley and Sons.

\bibitem[Ad\'{a}mek and Porst, 2003]{AP}
Ad\'{a}mek, J. and Porst, H. (2003).
\newblock On varieties and covarieties in a category.
\newblock {\em Mathematical Structures in Computer Science}, 13:201--232.

\bibitem[Borceux, 1994a]{bor1}
Borceux, F. (1994a).
\newblock {\em Encyclopedia of Mathematics an Its Applications: Handbook of
  Categorical Algebra 1.}
\newblock Cambridge University Press.

\bibitem[Borceux, 1994b]{bor2}
Borceux, F. (1994b).
\newblock {\em Encyclopedia of Mathematics an Its Applications: Handbook of
  Categorical Algebra 2.}
\newblock Cambridge University Press.

\bibitem[Kelly, 1980]{Kelly}
Kelly, G.~M. (1980).
\newblock A unified treatment of transfinite constructions for free algebras,
  free monoids, colimits, associated sheaves, and so on.
\newblock {\em Bulletin of the Australian Mathematical Society}, 22:1--83.

\bibitem[Lack, 2007]{lack}
Lack, S. (2007).
\newblock A 2-categories companion.
\newblock arXiv:math/0702535v1.

\bibitem[Lack and Street, 2002]{lackstr}
Lack, S. and Street, R. (2002).
\newblock The formal theory of monads 2.
\newblock {\em Journal of Pure and Applied Algebra}, 175:243--265.

\bibitem[Mac~Lane, 1971]{mac}
Mac~Lane, S. (1971).
\newblock {\em Categories for the Working Mathematician}.
\newblock Springer-Verlag.

\bibitem[Manes, 1976]{manes}
Manes, E.~G. (1976).
\newblock {\em Algebraic theories}.
\newblock Springer-Verlag.

\bibitem[Pavl\'{i}k, 2009]{p1}
Pavl\'{i}k, J. (2009).
\newblock Varieties defined without colimits.
\newblock In {\em Proceedings of the 7th Panhellenic Logic Symposium}.

\bibitem[Pavl\'{i}k, 2010a]{p}
Pavl\'{i}k, J. (2010a).
\newblock Free algebras in varieties.
\newblock {\em Archivum Mathematicum}, 46:25--38.

\bibitem[Pavl\'{i}k, 2010b]{pd}
Pavl\'{i}k, J. (2010b).
\newblock {\em On Categories of Algebras}.
\newblock PhD thesis, Masaryk University.

\bibitem[Rosick\'{y}, 1977]{Rac}
Rosick\'{y}, J. (1977).
\newblock On algebraic categories.
\newblock {\em Colloquia Mathematica Societatis J\'{a}nos Bolai}, 29:663--690.

\bibitem[Street, 1972]{street}
Street, R. (1972).
\newblock The formal theory of monads.
\newblock {\em Journal of Pure and Applied Algebra}, 2:149--168.

\end{thebibliography}

\end{document}